\RequirePackage{fix-cm} 
\documentclass[smallextended]{svjour3}       
\usepackage[utf8]{inputenc}
\smartqed                                                       
\DeclareMathAlphabet{\mathcal}{OMS}{cmsy}{m}{n}  
\usepackage{amsbsy,amsfonts,amsmath,amssymb,bbm,calc,caption,color,dsfont,graphics}
\usepackage{graphicx,ifthen,latexsym,mathptmx,mathrsfs}
\usepackage{overpic,pifont,psfrag,rotating,stmaryrd,theorem}
\usepackage{algorithm}
\usepackage{algorithmic}
\usepackage{bm}
\usepackage{float}
\usepackage{cite}
\usepackage{multirow}
\usepackage{multicol}
\usepackage[colorlinks,linkcolor=blue,citecolor=blue]{hyperref}
\usepackage{epstopdf}                   
\usepackage{geometry}                  
\usepackage{indentfirst}
\spnewtheorem{method}{Method}{\bf}{\it}

\usepackage[misc]{ifsym}

\numberwithin{equation}{section}                        
\numberwithin{theorem}{section}                        
\numberwithin{lemma}{section}                           
\numberwithin{corollary}{section}                       
\numberwithin{remark}{section}
\numberwithin{proposition}{section}
\numberwithin{problem}{section}

\begin{document}

\title{The Petty projection inequality for sets of finite perimeter}

\author{Youjiang Lin}

\authorrunning{Y. Lin} 

\institute{Y. Lin \at
             School of Mathematics and Statistics, Chongqing Technology and Business University, Chongqing  400067, PR China \\
              \email{yjl@ctbu.edu.cn}}


\maketitle

\begin{abstract}
The Petty projection inequality for sets of finite perimeter is proved. Our approach is based on Steiner symmetrization. Neither the affine
Sobolev inequality nor the functional Minkowski problem is used in our proof. Moreover, for sets of finite perimeter, we prove the Petty projection inequality with respect to Steiner symmetrization.
\end{abstract}
\footnotetext[1]{Research of the author is supported by National Natural Science Foundation of China NSFC 11971080 and Natural Science Foundation Projection of Chongqing cstc2018jcyjAX0790.}

\section{Introduction}
\label{intro}
\indent Within the Brunn-Minkowski theory, the two classical inequalities which connect the volume
of a convex body with that of its polar projection body are the Petty and Zhang projection inequalities. Unlike the classical isoperimetric inequality (see, e.g., \cite{MR3585546,MR0043486,MR2456887,MR2672283}), the Petty and Zhang projection inequalities are
affine isoperimetric inequalities in that they are inequalities between a pair of geometric functionals whose product is invariant under affine transformations. Many important results about affine isoperimetric inequalities and their functional forms have been found (see, e.g., \cite{MR3250365,MR3128983,MR3418448,MR2551138,MR2530600,MR2545028,MR4008006,MR2885586,MR1863023}).
The Petty projection inequality strengthens and directly implies the classical isoperimetric inequality, but it can be
viewed as an optimal isoperimetric inequality. It is the geometric core of the affine Sobolev-Zhang inequality \cite{MR1776095} which strengthens the classical sharp Euclidean Sobolev inequality. The $L_p$ version of Petty’s inequality by Lutwak et al. \cite{MR1863023} and its Orlicz extension by the same
authors \cite{MR2563216} both represent landmark results in the evolution of the Brunn-Minkowski
theory first towards an $L_p$ theory and, more recently, towards an Orlicz theory of convex
bodies. The $L_p$ and Orlicz theories of convex bodies have expanded rapidly (see, e.g.,  \cite{MR1901248,MR1898210,MR2251886,MR3194492,MR2652209,MR1725817,MR2680490,MR1942402,MR2159706,MR2067123,MR1378681,MR2785767,MR1987375,MR2211138,MR2652465,MR1231704,MR3488132,MR3128983,MR2652213,MR3695895,MR4018310,MR2873887}).

Affine isoperimetric inequalities referring to convex bodies have been extended to
certain classes of non-convex domains such as star bodies and sets of finite perimeter. The Petty projection inequality has been generalized first to compact domains with smooth boundary by Zhang \cite{MR1776095} and  to sets of finite perimeter by Wang \cite{MR2927377}. Recently, Haberl and Schuster \cite{MR4008006} show every even, zonal measure on the Euclidean
unit sphere gives rise to an isoperimetric inequality for
sets of finite perimeter which directly implies the classical
Euclidean isoperimetric inequality. The strongest member of
this large family of inequalities was shown to be the only affine
invariant one among them--the Petty projection inequality. In \cite{MR2927377}, Wang used the functional
Minkowski problem on $BV(\mathbb{R}^n)$ and an approach to the affine Sobolev-Zhang inequality for compactly supported $C^1$ functions developed by Lutwak et al. \cite{MR2211138} to
establish the affine Sobolev-Zhang inequality of $BV$ functions. As a consequence, the Petty projection inequality for sets of finite perimeter was established. Using the same ideas, Haberl and Schuster \cite{MR4008006} obtained the large family of sharp Sobolev type
inequalities which can be seen as the functional form of the large family of isoperimetric inequalities.

In this paper, we give a direct proof of the Petty projection inequality for sets of finite perimeter. Our
approach is based on Steiner symmetrization. Neither the affine Sobolev inequality nor the functional Minkowski problem is used in our proof.  Moreover, for sets of finite perimeter, we prove the Petty projection inequality with respect to Steiner symmetrization.

Let $\mathcal{C}^n$ denote the set of  compact sets of finite perimeter in $\mathbb{R}^n$ with nonempty interior. For $E\in\mathcal{C}^n$, let $\partial^{\ast}E$ and $\nu^E=(\nu_1^E,\dots,\nu_{n-1}^E,\nu_y^E)$ denote the reduced boundary and the generalized inner normal to $E$, respectively (see Section 3.2 for precise
definitions). We define the projection body $\Pi E$ of $E$ to be the convex set with support function
\begin{eqnarray}\label{1e}
h(\Pi E,z)=\frac{1}{2}\int_{\partial^{\ast}E} |z\cdot \nu^E(x)| d\mathcal{H}^{n-1}(x),\;\;z\in\mathbb{R}^n.
\end{eqnarray}
This extension of Minkowski's classical notion of the projection body of a convex body
was first given by Wang \cite{MR2927377}, who generalized a definition of Zhang \cite{MR1776095} for compact sets
with piecewise $C^1$ boundary.

For the polar (see Section \ref{SS1} for definitions) of $\Pi E$ we will write $\Pi^{\ast}E$. The $\mathcal{L}^n$ Lebesgue measure of a set of finite perimeter $E\in\mathcal{C}^n$ will be denoted by $|E|$, and for the volume of the unit ball in $\mathbb{R}^n$ we use $\omega_n$.

We assume that $E\in\mathcal{C}^n$ satisfies the condition
\begin{eqnarray}\label{1b}
\mathcal{H}^{n-1}(\{x\in\partial^{\ast}E:\;\nu_y^{E}(x)=0\})=0.
\end{eqnarray}

\begin{theorem}\label{T1}
Let $E\in\mathcal{C}^n$ satisfy (\ref{1b}) and $E^s$ denote its Steiner symmetrization. Then
\begin{eqnarray}\label{1a}
|E|^{n-1}|\Pi^{\ast}E|\leq |E^s|^{n-1}|\Pi^{\ast}E^s|.
\end{eqnarray}
\end{theorem}

\begin{theorem}\label{T2}
Let $E\in\mathcal{C}^n$. Then
\begin{eqnarray}\label{1i}
|E|^{n-1}|\Pi^{\ast}E|\leq(\omega_n/\omega_{n-1})^n.
\end{eqnarray}
\end{theorem}

Note that since the absolute value function $|\cdot|$ is not strictly convex, we cannot characterize the cases of equality in (\ref{1a}) and (\ref{1i}) using the methods in the paper.

In section \ref{S2}, we set up notation and terminology and compile some basic facts about the Brunn-Minkowski theory of convex bodies and the theory of functions of bounded variation. In Section \ref{S3}, we prove that for any a set of finite perimeter $E$, there exists
a sequence of unit directions $\{u_i\}$ such that there exists a subsequence of the successive Steiner symmetrizations $E_i:=S_{u_i}\cdots S_{u_1}E$ of $E$ which converges to an origin-symmetric ball with same volume. The main difficulty in proving the convergence is that every direction $u_i$ of Steiner symmetrization is related to $E_{i-1}$. Section \ref{S4} is devoted to the study of some basic properties of projection bodies of sets of finite perimeter. For example, the continuity and the affine invariance of the Petty projection operator. It is worth pointing out that the proofs of these properties are new and  did not appear in the known literatures.

In Section \ref{S5}, we give the proofs of main theorems. First, we prove the monotonicity of the volumes of projection bodies with respect to Steiner symmetrization. To prove the monotonicity, we make critical
use of a beautiful fact found by Lutwak et al. in \cite{MR2563216}. The fact gives a useful way to imply the inclusion relation between the Steiner symmetrization of a convex body and another convex body. In the remarkable paper \cite{MR2178968}, Chleb\'ik et al. characterized the sets of finite perimeter whose perimeter
is preserved under Steiner symmetrization and proved the monotonicity of the
perimeter of sets of finite perimeter with respect to Steiner symmetrization. The ideas and techniques of Chleb\'ik et al. \cite{MR2178968} play a critical role in this paper. Comparing to the paper \cite{MR2178968}, the perimeter of sets of finite perimeter is replaced by the volume of the projection body of sets of finite perimeter, which is a function of the volume of the unit ball of an $n$-dimensional Banach space. Finally, using the convergence of Steiner symmetrization, monotonicity and continuity of the projection operator, we prove the Petty projection inequality for sets of finite perimeter.
\section{Preliminaries}
\label{S2}
Let $\mathbb{N}$ denote the set of positive integers. Let $\{e_1$,\dots,$e_n\}$ denote the standard orthonormal basis of the Euclidean space $\mathbb{R}^n$.  A point $x\in\mathbb{R}^n$, $n\geq 2$, will be usually labeled by $(x^{\prime},y)$, where $x^{\prime}=(x_1,\dots,x_{n-1})\in\mathbb{R}^{n-1}$ and $y\in\mathbb{R}$. To emphasize the different role of the variable $y$, we shall also write $\mathbb{R}^n=\mathbb{R}^{n-1}\times\mathbb{R}_y$. We shall use $\mathbb{R}^+$ and $\mathbb{R}^-$ to denote $[0,+\infty)$ and $(-\infty,0]$, respectively. For $u\in S^{n-1}$, let $u^{\perp}$  denote the codimension $1$ subspace of $\mathbb{R}^n$
that is orthogonal to $u$. If $E$ is a measurable subset of $\mathbb{R}^n$ and $E$ is contained in an $i$-dimensional affine subspace of $\mathbb{R}^n$ but in no affine subspace of lower dimension, then $|E|$ will denote the $i$-dimensional Lebesgue measure of $E$. Two measurable sets $E$ and $F$ are equivalent if the Lebesgue measure of their symmetric
 difference $E \triangle F$ is zero. Let $B^n$ denote the Euclidean unit ball centered at the origin in $\mathbb{R}^n$. Let $B_r(x)$ denote the ball, centered at $x$, having radius $r$. If $x\in\mathbb{R}^n$ then by abuse of notation we will write
$|x|=\sqrt{x\cdot x}$.
\subsection{Convex bodies.}
\label{SS1}
We develop some notation and, for quick later reference, list some basic facts about convex bodies. Good general references for the theory of convex bodies are provided by the books of Gardner \cite{MR2251886}, Gruber \cite{MR1242983,MR2335496}, Webster \cite{MR1443208} and Schneider \cite{MR3155183}.

 We write $\mathcal{K}^n$ for the set of convex bodies (compact convex subsets) of $\mathbb{R}^n$. We write $\mathcal{K}^n_0$ for the set of convex bodies that contain the origin in their interiors. For $K\in \mathcal{K}^n$, let $h(K;\cdot)=h_K:\mathbb{R}^n\rightarrow\mathbb{R}$ denote the {\it support function} of $K$; i.e.,
$$h(K;x)=\max\{x\cdot z:\;z\in K\}.$$

For $K\in \mathcal {K}_0^n$, its {\it polar body} $K^{\ast}$ is defined by
$$K^{\ast}=\{x\in\mathbb{R}^n:x\cdot z\leq 1\;{\rm for\;all}\;z\in K\}.$$

For $A\in SL(n)$ and $K\in\mathcal{K}_0^n$,
\begin{eqnarray}\label{1f}
(AK)^{\ast}=A^{-t}K^{\ast}.
\end{eqnarray}

For $K\in \mathcal{K}_0^n$, its {\it radial function} $\rho_K$ is defined as
\begin{eqnarray}\label{2kk}
\rho_K(x)=\max\{\lambda\geq 0:\;\lambda x\in K\},\;\;x\in\mathbb{R}^n\backslash\{0\}.
\end{eqnarray}
It is easily verified that
\begin{eqnarray}\label{2k}
\rho_{K^{\ast}}(x)=\frac{1}{h_K(x)},\;\;x\in\mathbb{R}^n\setminus\{0\}.
\end{eqnarray}
\subsection{Functions of bounded variation and Sets of finite perimeter}
In this section, we review some basic definitions and facts about functions of bounded variation and sets of finite perimeter on $\mathbb{R}^n$.  Good general references are Ambrosio et al. \cite{MR1857292}, Cianchi \cite{MR2723814}, Evans \cite{MR3409135}, Ma${\rm z}^{\prime}$ya \cite{MR817985},\cite{MR2777530} and Ziemer \cite{MR1014685}.

 The space of functions of bounded variation in $\mathbb{R}^n$ is denoted by $BV(\mathbb{R}^n)$. Recall that a function $f\in L^1(\mathbb{R}^n)$ is said to be of bounded variation in $\mathbb{R}^n$ if it is integrable  and its distributional gradient $Df$ is a vector-valued Radon measure in $\mathbb{R}^n$ whose total variation $|Df|$ is finite in $\mathbb{R}^n$;
thus
$$BV(\mathbb{R}^n)=L^1(\mathbb{R}^n)\cap\{f:\;Df\;{\rm is\;a\;measure},\;|Df|(\mathbb{R}^n)<\infty\}.$$

For $f\in  L^1(\mathbb{R}^n)$, let $f^{\ast}$ denote the {\it precise representative} of $f$, i.e.

\begin{equation}
f^{\ast}(x)=
\left\{ \begin{aligned}
         &\lim_{r\rightarrow 0}\frac{1}{r^n\omega_n}\int_{B_r(x)}f(z)dz&{\rm if\;this\;limit\;exists} \\
                 & 0&{\rm otherwise}.
                          \end{aligned} \right.
                          \end{equation}

For $f_i,f\in BV(\mathbb{R}^n)$, $i\in\mathbb{N}$, we say that {\it $f_i$ weakly$^{\ast}$ converges in $BV(\mathbb{R}^n)$ to $f$} if $f_i$ converges to $f$ in $L^1(\mathbb{R}^n)$ and $Df_i$ weakly$^{\ast}$ converges to $Df$ in $\mathbb{R}^n$.

A measurable subset $E$ of $\mathbb{R}^n$ is said to be of {\it finite perimeter} if its characteristic function $\chi_E$ is a bounded variation function in $\mathbb{R}^n$. Let $\mathcal{C}^n$ denote the set of compact sets of finite perimeter in $\mathbb{R}^n$ with nonempty interior. The perimeter of $E$ in $\mathbb{R}^n$ denoted by $P(E)$ is defined by
\begin{eqnarray}\label{2n}
P(E)=|D\chi_E|(\mathbb{R}^n).
\end{eqnarray}

The {\it Hausdorff distance} of the sets $E,F\in\mathcal{C}^n$ is defined by
$$d_H(E,F)=\min\{\lambda\geq 0:\;K\subset L+\lambda B^n,\;L\subset K+\lambda B^n\}.$$

For $K,L\in\mathcal{K}^n$, the  Hausdorff distance of $K$ and $L$ (see, e.g., \cite[Lemma 1.8.14]{MR3155183})
$$d_H(K,L)=\sup_{u\in S^{n-1}}\left|h_K(u)-h_L(u)\right|.$$

For $E,F\in\mathcal{C}^n$, we define $L_1$ or also the {\it symmetric-difference distance} (sometimes also called Nikod\'ym distance, as in \cite{MR1242983}) as follows,
$$d_1(E,F)=\mathcal{L}^n(E\triangle F).$$

Let $E\in\mathcal{C}^n$ and let $D_i\chi_E$ denote the $i$-th component of the distributional gradient $D\chi_E$. We denote by $\nu_i^E$, $i=1,\dots,n$, the derivative of the measure $D_i\chi_E$ with respect to $|D\chi_E|$. The {\it reduced boundary} $\partial^{\ast}E$ of $E$ is the set of all points $x\in \mathbb{R}^n$ such that the vector $\nu^E(x)=(\nu_1^E(x),\dots,\nu^E_{n}(x))$ exists and satisfies $|\nu^E(x)|=1$. The vector $\nu^E(x)$ is called the {\it generalized inner normal} to $E$ at $x$.  Throughout this paper, let $\nu^E_{x'}$ and $\nu^E_y$ denote the front $n-1$ components and the $n$-th component of $\nu^E$, respectively.

We have (see \cite[Theorem 3.59]{MR1857292})
\begin{eqnarray}\label{2b}
D\chi_E=\nu^E\mathcal{H}^{n-1}\llcorner\partial^{\ast}E.
\end{eqnarray}
Eq. (\ref{2b}) implies that
\begin{eqnarray}\label{2c}
|D\chi_E|=\mathcal{H}^{n-1}\llcorner\partial^{\ast}E
\end{eqnarray}
and that
$$|D_i\chi_E|=|\nu_i^E|\mathcal{H}^{n-1}\llcorner\partial^{\ast}E,\;\;i=1,\dots,n.$$

For any $E\in\mathcal{C}^n$ and $u\in S^{n-1}$, define, for $x^{\prime}\in u^{\perp}$,
$$E_{x',u}=\{x'+tu:\;t\in\mathbb{R},\;x'+tu\in E\}$$
and
$$\ell_{E,u}(x^{\prime})=|E_{x^{\prime},u}|.$$
For $u=e_n$, let $E_{x'}$ and $\ell_E$ denote $E_{x^{\prime},e_n}$ and $\ell_{E,e_n}$, respectively. $E$ will be suppressed when clear from the context, and thus we
will often denote $\ell_E(x')$ by $\ell(x')$.

Let $\pi_u(E)$ denote the orthogonal projection of $E\subset\mathbb{R}^{n}$ onto $u^{\perp}$.
 In what follows, the essential projection of a set $E\subset\mathbb{R}^{n}$ onto $u^{\perp}$ is defined as
$$\pi_u(E)^+=\{x^{\prime}\in u^{\perp}:\;\ell_{E,u}(x^{\prime})>0\}.$$
For $u=e_n$, let $\pi_{n-1}(E)$ and $\pi_{n-1}(E)^+$ denote $\pi_{e_n}(E)$ and $\pi_{e_n}(E)^+$, respectively.

The following theorem is a special case of the co-area formula for rectifiable sets (see\cite[Theorem 2.93 and Remark 2.94]{MR1857292}).

\begin{theorem}\label{T8}
Let $E\in\mathcal{C}^n$. Let $g:\mathbb{R}^{n-1}\rightarrow [0,+\infty)$ be any Borel function. Then
\begin{eqnarray}\label{2l}
\int_{\partial^{\ast}E}g(x)|\nu^E_y(x)|d\mathcal{H}^{n-1}(x)=\int_{\mathbb{R}^{n-1}}dx^{\prime}\int_{(\partial^{\ast}E)_{x^{\prime}}}g(x^{\prime},y)d\mathcal{H}^0(y).
\end{eqnarray}
\end{theorem}

Next, we give a theorem concerning one-dimensional sections of sets of finite perimeter, it can be easily deduced from \cite[Theorem 3.108]{MR1857292}.

\begin{theorem}\label{T7}
Let $E$ be a set of finite perimeter in $\mathbb{R}^n$. Then, for $\mathcal{L}^{n-1}$-a.e. $x^{\prime}\in\mathbb{R}^{n-1}$,
\begin{eqnarray}\label{2e}
E_{x^{\prime}}\; has\;finite\;perimeter\;in\;\mathbb{R}_y,
\end{eqnarray}
\begin{eqnarray}\label{2ee}
(\partial^{\ast}E)_{x^{\prime}}=\partial^{\ast}(E_{x^{\prime}}),
\end{eqnarray}
\begin{eqnarray}\label{2gg}
\nu^E_y(x^{\prime},y)\neq 0\; for\;every\;y\; such\;that\;(x^{\prime},y)\in\partial^{\ast}E,
\end{eqnarray}
\begin{equation}\label{2f}
\left\{ \begin{aligned}
\lim_{\eta\rightarrow y^+}\chi^{\ast}_E(x^{\prime},\eta)=1,\;\;\lim_{\eta\rightarrow y^-}\chi^{\ast}_E(x^{\prime},\eta)=0\;\;{\rm if}\;\nu_y^E(x^{\prime},y)>0,\\
\lim_{\eta\rightarrow y^+}\chi^{\ast}_E(x^{\prime},\eta)=0,\;\;\lim_{\eta\rightarrow y^-}\chi^{\ast}_E(x^{\prime},\eta)=1\;\;{\rm if}\;\nu_y^E(x^{\prime},y)<0.
\end{aligned} \right.
\end{equation}
In particular, a Borel set $G_E\subset\pi_{n-1}(E)^+$ exists such that $\mathcal{L}^{n-1}(\pi_{n-1}(E)^+\setminus G_E)=0$ and  (\ref{2e})-(\ref{2f}) are fulfilled for every $x^{\prime}\in G_E$.
\end{theorem}

\begin{remark}\label{3s}
By \cite[Prop. 3.52]{MR1857292}, if $E\subset \mathbb{R}$ is a measurable set, then $E$ has finite perimeter in $\mathbb{R}$ if and only if there exist $-\infty\leq a_1<b_1<\dots<b_m\leq+\infty$ such that
$$E=\bigcup_{i=1}^{m}(a_i,b_i)$$
up to a set of measure zero.

In Theorem \ref{T7}, if $E$ is a set of finite perimeter in $\mathbb{R}^n$ and $\mathcal{\ell}_E(x')<\infty$ for $\mathcal{L}^{n-1}$-a.e. $x'\in\mathbb{R}^{n-1}$, then for every $x^{\prime}\in G_E$,  $E_{x^{\prime}}$ has finite perimeter in $\mathbb{R}^n_{x^{\prime}}$ and $E_{x^{\prime}}$ is bounded. Thus,
$$E_{x^{\prime}}\;\; is\; equivalent\; to\;\{x'\}\times\bigcup^{m(x^{\prime})}_{i=1}(a_i,b_i)$$
where
$$m(x^{\prime})=\frac{\mathcal{H}^0(\partial^{\ast}(E_{x^{\prime}}))}{2},\;\;\;\;(x',a_i),\;(x',b_i)\in\partial^{\ast}(E_{x^{\prime}})$$
and
$$\nu_y^E(x^{\prime},a_i)>0\;\;and\;\;\nu_y^E(x^{\prime},b_i)<0.$$
\end{remark}

The approximation theorem was proved by Maggi in \cite[Theorem 2.4]{MR2402947}.

\begin{theorem}(Density of smooth or polyhedral sets).\label{T11}
If $E$ is a set of finite perimeter in $\mathbb{R}^n$, then a
sequence $E_i$ of bounded open sets with with smooth or polyhedral boundary can be found so that $E_i\rightarrow E$
in the $L^1$ distance and $P(E_i)\rightarrow P(E)$.
\end{theorem}
\subsection{Steiner symmetrization and spherical symmetrization}
Steiner symmetrization is a classical and very well-known
device, which has seen a number of remarkable applications to problems of
geometric and functional nature (see, e.g., \cite{MR3581298,MR1475547,MR1394967,MR2079885,MR2178968,MR1947097,MR2228056,MR3723145,MR1322313,MR1465371,MR3018162}).

 For $E\in\mathcal{C}^n$ and $u\in S^{n-1}$, we define the {\it Steiner symmetrization}, $S_uE$, of $E$ about the hyperplane $u^{\perp}$ as
$$S_uE=\{x^{\prime}+tu\in\mathbb{R}^n:x'\in\pi_u(E),\;|t|\leq\ell_{E,u}(x^{\prime})/2\}.$$
For $u=e_n$, let $SE$ (or $E^s$) denote $S_{e_n}E$.

For $E\in\mathcal{C}^n$, we define the {\it spherical symmetrization}, $E^{\star}$, of $E$ is the close ball centered at the origin which has the same Lebesgue measure as $E$.

Specially, for $K\in\mathcal {K}_0^n$, we have
$$K^s=\left\{\left(x^{\prime},\frac{1}{2}y_1+\frac{1}{2}y_2\right)\in\mathbb{R}^{n-1}\times\mathbb{R}_y:(x^{\prime},y_1),\;(x^{\prime},-y_2)\in K\right\}.$$

In this paper, we shall make critical use of the following fact that was provided by Lutwak et al. in \cite[Lemma 1.1]{MR2563216}.
\begin{lemma}\label{L2}
Suppose $K,L\in\mathcal {K}_0^n$ and consider $K,L\subset\mathbb{R}^{n-1}\times\mathbb{R}_y$. Then
$$SK^{\ast}\subset L^{\ast},$$
if and only if
$$h_K(x^{\prime},t)=1=h_K(x^{\prime},-s),\;\; with\;t\neq-s\;\;\Longrightarrow \;\;h_L(x^{\prime},t/2+s/2)\leq 1.$$
In addition, if $SK^{\ast}=L^{\ast}$, then $h_K(x^{\prime},t)=1=h_K(x^{\prime},-s)$, with $t\neq-s$ implies that $h_L(x^{\prime},t/2+s/2)= 1$.
\end{lemma}

The following properties of Steiner symmetrization were summarized by Talenti in \cite[P.107]{MR1242977}.

\begin{theorem}\label{T9}
If  $E\in\mathcal{C}^n$, then $E^s\in\mathcal{C}^n$ and
$$P(E^s)\leq P(E).$$
\end{theorem}

\section{The convergence of Steiner symmetrization}\label{S3}
For $E\in\mathcal{C}^n$ and $u\in S^{n-1}$. Let
$$\partial^{\ast}_u E=\{x\in\partial^{\ast}E:\;u\cdot\nu^{E}(x)=0\}$$
and
\begin{eqnarray}\label{3d}
T(E)=\{u\in S^{n-1}:\;\mathcal{H}^{n-1}(\partial_u^{\ast}E)=0\}.
\end{eqnarray}

\begin{theorem}\label{T3}
Let $E\in \mathcal{C}^n$. Then there exits $\{u_i\}_{i=1}^{\infty}\subset S^{n-1}$ such that $u_1\in T(E)$,
$u_i\in T(E_{i-1})$ when $i\geq2$, where $E_i=S_{u_i}\dots S_{u_1}E$, and there exists a subsequence of $\{E_i\}$, denoted by $\{E_{i_j}\}$, satisfying $E_{i_j}\rightarrow E^{\star}$
in the Hausdorff distance.
\end{theorem}

In order to prove Theorem \ref{T3}, we need the following Lemmas \ref{L12}-\ref{L13}.

\begin{lemma}\cite[Lemma 6.12]{MR4018310}\label{L12}
Let $E\in\mathcal{C}^n$ and $T(E)$ as in (\ref{3d}).
Then
$$\mathcal{H}^{n-1}\left(S^{n-1}\backslash T(E)\right)=0.$$
\end{lemma}

\begin{lemma}\label{L6}
Let $E_i\in\mathcal{C}^n$, $i\in\mathbb{N}$, and let $E$ be a compact set in $\mathbb{R}^n$. If
$$d_H(E_i,E)\rightarrow 0\;\;when\;\;i\rightarrow\infty,$$
then
$$\max_{x\in E_i}|x|\rightarrow\max_{x\in E}|x|.$$
\end{lemma}
\noindent{\bf Proof.} Since $d_H(E_i,E)\rightarrow 0$ when $i\rightarrow\infty$, for any $\varepsilon>0$, there exists a positive
integer $N_{\varepsilon}$ such that for any $i>N_{\varepsilon}$,
$$E_i\subset E+\varepsilon B^n\;\;{\rm and}\;\;E\subset E_i+\varepsilon B^n.$$
Thus
$$\max_{x\in E_i}|x|\leq\max_{x\in E+\varepsilon B^n}|x|\leq \max_{x\in E}|x|+\varepsilon$$
and
$$\max_{x\in E}|x|\leq\max_{x\in E_i+\varepsilon B^n}|x|\leq \max_{x\in E_i}|x|+\varepsilon.$$
Therefore,
$$\left|\max_{x\in E}|x|-\max_{x\in E_i}|x|\right|\leq \varepsilon.$$
This completes the proof of the lemma.\qed

\begin{corollary}\label{C1}
Let $E$, $E_i$, $i\in\mathbb{N}$, be compact sets satisfying $E_{i+1}\subset E_i$ and
\begin{eqnarray}\label{2d}
E=\bigcap_{i=1}^{\infty}E_i.
\end{eqnarray} Then
$$\max_{x\in E_i}|x|\rightarrow\max_{x\in E}|x|.$$
\end{corollary}
\noindent{\bf Proof.}
By (\ref{2d}) and the compactness of $E_i$, the sequence $E_i$ converges to $E$ in the Hausdorff distance. This and Lemma \ref{L6} yield the desired result. \qed

\begin{lemma}\label{L13}
Let $E$, $E_i$, $i\in\mathbb{N}$, be compact sets satisfying $E_{i+1}\subset E_i$ and
$$E=\bigcap_{i=1}^{\infty}E_i.$$ Then
\begin{eqnarray}\label{1l}
SE=\bigcap_{i=1}^{\infty}SE_i.
\end{eqnarray}
\end{lemma}
\noindent{\bf Proof.}
First, we prove that the orthogonal projections of $SE$ and $\cap_{i=1}^{\infty}SE_i$ onto $\mathbb{R}^{n-1}$ are same, i.e.,
\begin{eqnarray}\label{1h}
\pi_{n-1}(SE)=\pi_{n-1}(\cap_{i=1}^{\infty}SE_i).
\end{eqnarray}
If $x'\in \pi_{n-1}(\cap_{i=1}^{\infty}SE_i)$, then for any $i$, $x'\in \pi_{n-1}(SE_i)$. Since $\pi_{n-1}(E_i)=\pi_{n-1}(SE_i)$, $x'\in \pi_{n-1}(E_i)$. Thus there exists $x^i\in E_i$, $i=1,2,\dots$, such that the orthogonal projection of $x^i$ onto $\mathbb{R}^{n-1}$ is $x'$. Since $E_{i+1}\subset E_i$, $x^i\subset E_1$ for any $i\in\mathbb{N}$. Thus for the sequence $\{x^i\}$, there exist a subsequence $\{x^{i_j}\}$ and  $x^0\in\mathbb{R}^n$ such that
$$x^{i_j}\rightarrow x^0\;\;{\rm when}\;j\rightarrow \infty.$$
It is clear that the orthogonal projection of $x^0$ onto $\mathbb{R}^{n-1}$ is also $x'$.
Since $E_{i+1}\subset E_i$ and the compactness of $E_i$, we have $x^0\in E_i$ for any $i$. Thus, $x^0\in \cap_{i=1}^{\infty}E_i$. This implies
$$x'\in \pi_{n-1}(\cap_{i=1}^{\infty}E_i)=\pi_{n-1}(E)=\pi_{n-1}(SE).$$ Therefore,
$$\pi_{n-1}(\cap_{i=1}^{\infty}SE_i)\subset\pi_{n-1}(SE).$$
Moreover, $\pi_{n-1}(SE)\subset\pi_{n-1}(\cap_{i=1}^{\infty}SE_i)$ is clear. Thus, (\ref{1h}) is established.

 Next, we prove that for any $x'\in \pi_{n-1}(SE)$,
\begin{eqnarray}\label{1j}
(SE)_{x'}=\left(\cap_{i=1}^{\infty}SE_i\right)_{x'}.
\end{eqnarray}
Since $E_{i+1}\subset E_i$, $(E_{i+1})_{x'}\subset (E_i)_{x'}$. By the limit theorem with respect to  sequences of measurable sets (see \cite[Theorem 1.2 (iv)]{MR3409135}), we have
\begin{eqnarray}\label{1k}
|(SE)_{x'}|=|(E)_{x'}|=|(\cap_{i=1}^{\infty}E_i)_{x'}|=|\cap_{i=1}^{\infty}(E_i)_{x'}|=\lim_{i\rightarrow\infty}|(E_i)_{x'}|=\lim_{i\rightarrow\infty}|(SE_i)_{x'}|=|\left(\cap_{i=1}^{\infty}SE_i\right)_{x'}|.\nonumber
\end{eqnarray}
This and the symmetry and compactness of $SE$  and $\cap_{i=1}^{\infty}SE_i$ yield (\ref{1j}). The desired equality (\ref{1l}) now follows from (\ref{1h}) and (\ref{1j}).\qed
\

\noindent{\bf Proof of Theorem \ref{T3}.}
Define
$$\Gamma_E=\{S_{u_k}\dots S_{u_1}E:\;k\in\mathbb{N},\;u_1\in T(E),\;u_2\in T(E_1),\dots,u_k\in T(E_{k-1})\},$$
where $T(E)$ is defined as in (\ref{3d}) and $E_j=S_{u_j}\dots S_{u_1}E$, $j=1,\dots,k-1$.

For $E\in \mathcal{C}^n$, let $r_E=\max_{x\in E}|x|$, which is the minimal radius of
balls centered at the origin that contain $E$. Let $r_1$ be the infimum of all $r_C$, where $C\in \Gamma_E$. Then there is a sequence of $\{C_i\}\subset\Gamma_E$ so that $r_{C_i}\rightarrow r_1$. Obviously, the sequence $\{C_i\}$ is bounded, because each $C_i\subset r_EB^n$. By \cite[Theorem 1.8.5]{MR3155183}, there is a subsequence $C_{i_k}$ that converges to a compact set $\bar{E}$ in the Hausdorff distance.  By Lemma \ref{L6}, $r_{\bar{E}}=r_1$. Denote $r_1B^n$ by $B_1$, it is clear that $\bar{E}\subset B_1$.

Next we prove $\bar{E}=B_1$. Assume it is not true. There is a
small open cap $U$ on $\partial B_1$ so that $U\cap \bar{E}=\emptyset$. For any line
$\xi$ such that $\xi\cap U\neq \emptyset$, either $\xi\cap \bar{E}=\emptyset$ or the
line $\xi$ intersects a longer chord in $B_1$ than in $\bar{E}$; that is,
$|\xi\cap B_1|>|\xi\cap \bar{E}|$. After taking a Steiner symmetrization $S_u\bar{E}$
for some $u\in S^{n-1}$, the symmetrization $S_u\bar{E}$ fails to intersect both $U$
and a new cap $U^{\prime}$ given by the reflection of $U$ with respect to
the hyperplane $u^{\perp}$. One can continue to take symmetrizations with respect to an appropriate finite
family of hyperplanes with normals $v_1,\dots,v_m\in S^{n-1}$
that generate finitely many caps covering the whole sphere $\partial B_1$ and
generate a compact set of finite perimeter about origin
$\tilde{E}=S_{v_m}\dots S_{v_1}\bar{E}$ so that $|\xi\cap B_1|>|\xi\cap \tilde{E}|$ for any line
such that $\xi\cap\partial B_1\neq\emptyset$. Thus, $r_{\tilde{E}}<r_1$.

By the above analysis,
$S^{n-1}\subset \bigcup_{i=1}^{m} U_{v_i}$,
where $U_{v_i}$ is the reflection of $U$ with respect to the hyperplane $v_i^{\perp}$. Since the cap $U$ is open and $S^{n-1}$ is compact, there exist sufficiently mall $\delta_1,\dots,\delta_m>0$ such that any $u_i\in B_{\delta_i}(v_i)\cap S^{n-1}$, $i=1,\dots,m$, we have
$S^{n-1}\subset \bigcup_{i=1}^{m} U_{u_i}$. By the analysis of the above paragraph, for any  $u_i\in B_{\delta_i}(v_i)\cap S^{n-1}$, $i=1,\dots,m$, let $\bar{E}_1=S_{u_m}\dots S_{u_1}\bar{E}$. Then
\begin{eqnarray}\label{3ll}
r_{\bar{E}_1}<r_1.
\end{eqnarray}

For the above $\delta_1,\delta_2,\dots,\delta_m$, by Lemma \ref{L12}, there exits $u_1\in B_{\delta_1}(v_1)\cap S^{n-1}$ such that
$$u_1\in \bigcap_{k=1}^{\infty}T(C_{i_k}).$$
Similarly, there exists $u_2\in  B_{\delta_2}(v_2)\cap S^{n-1}$ such that
$$u_2\in \bigcap_{k=1}^{\infty}T(S_{u_1}C_{i_k}).$$
Continue to take the process, we can get $u_1,\dots,u_m$ such that $u_j\in B_{\delta_j}(v_j)\cap S^{n-1}$, $j=1,\dots,m$, and
$$u_j\in \bigcap_{k=1}^{\infty}T(S_{u_{j-1}}\dots S_{u_1}C_{i_k}).$$

Denote $\tilde{C}_{i_k}=S_{u_m}\cdots S_{u_1}C_{i_k}$. Since
$C_{i_k}\rightarrow \bar{E}$ in the Hausdorff distance. For any positive integer $i$, there exists a positive integer $N$ such that any $k>N$
$$C_{i_k}\subset \bar{E}+\frac{1}{i}B^n.$$
Thus
$$S_{u_m}S_{u_{m-1}}\dots S_{u_1}C_{i_k}\subset S_{u_m}S_{u_{m-1}}\dots S_{u_1}\left(\bar{E}+\frac{1}{i}B^n\right),$$
which implies that
\begin{eqnarray}\label{3l}
\max\left\{|x|:\;x\in S_{u_m}S_{u_{m-1}}\dots S_{u_1}\left(\bar{E}+\frac{1}{i}B^n\right)\right\}\geq \max\{|x|:\;x\in S_{u_m}S_{u_{m-1}}\dots S_{u_1}C_{i_k}\}\geq r_1.
\end{eqnarray}
By Lemma \ref{L13}, Corollary \ref{C1} and (\ref{3l}), we have
\begin{eqnarray*}
&&\max\{|x|:\;x\in S_{u_m}S_{u_{m-1}}\dots S_{u_1}\bar{E}\}\nonumber\\
&=&\max\left\{|x|:\;x\in S_{u_m}S_{u_{m-1}}\dots S_{u_1}\bigcap_{i=1}^{\infty}\left(\bar{E}+\frac{1}{i}B^n\right)\right\}\nonumber\\
&=&\max\left\{|x|:\;x\in \bigcap_{i=1}^{\infty}S_{u_m}S_{u_{m-1}}\dots S_{u_1}\left(\bar{E}+\frac{1}{i}B^n\right)\right\}\nonumber\\
&=&\lim_{i\rightarrow\infty}\max\left\{|x|:\;x\in S_{u_m}S_{u_{m-1}}\dots S_{u_1}\left(\bar{E}+\frac{1}{i}B^n\right)\right\}\nonumber\\
&\geq&r_1.
\end{eqnarray*}
Let $\bar{E}_1=S_{u_m}\dots S_{u_1}\bar{E}$. Then
$r_{\bar{E}_1}\geq r_1$. This contradicts (\ref{3ll}).

We have shown that for any $E\in \mathcal{C}^n$, there are
$S_{u_{i_1}}\cdots S_{u_1}E\in\Gamma_E$ so that the Hausdorff distance between
$S_{u_{i_1}}\cdots S_{u_1}E$ and the centered ball $B_1$ can be arbitrarily small.

For a sequence of positive numbers $\varepsilon_k\rightarrow 0^+$, there is
$D_1:=S_{u_{i_1}}\cdots S_{u_1}E\in\Gamma_E$ so that $d_H(D_1,B_1)<\varepsilon_1$. Similarly, there are
$D_2:=S_{u_{i_2}}\dots S_{u_{i_1+1}} D_1\in\Gamma_{D_1}$ so that $d_H(D_2,B_1)<\varepsilon_2$. In general, for $k=3,4,\dots$, there are
$D_k:= S_{u_{i_k}}\dots S_{u_{i_{k-1}+1}}D_{k-1}\in\Gamma_{D_{k-1}}$ so that $d_H(D_k,B_1)<\varepsilon_k$.
Continue the  process, we can get a sequence $\{D_k\}_{k=1}^{\infty}$ and
$D_k\rightarrow B_1$ in the Hausdorff distance.
Let
$$E_j=S_{u_j}S_{u_{j-1}}\dots S_{u_1}E,\;\;j\in\mathbb{N}.$$
Then $\{D_k\}$ is a subsequence of $\{E_j\}$ while $\{D_k\}$ and $\{E_j\}$ satisfy the conclusions of the theorem. \qed

\section{Definition and basic properties of projection bodies}\label{S4}

For $E\in\mathcal{C}^n$, we define the projection body $\Pi E$ of $E$ to be the convex set with support function
\begin{eqnarray}\label{3n}
h_{\Pi E}(z)=\frac{1}{2}\int_{\partial^{\ast}E} |z\cdot \nu^E(x)| d\mathcal{H}^{n-1}(x),\;\;z\in \mathbb{R}^n.
\end{eqnarray}
It is clear that the function $h_{\Pi K}(\cdot)$
is the support function of a convex body, $\Pi K$, that contains the origin in
its interior. The polar body of $\Pi E$ will be denoted by $\Pi^{\ast}E$, rather than $(\Pi E)^{\ast}$.

The following proposition shows that the Petty projection operator $\Pi:\mathcal{C}^n\rightarrow\mathcal{K}_0^n$ is
continuous when $E_i$ converges to $E_0$ in the $L_1$ distance and $P(E_i)\rightarrow P(E_0)$.

\begin{proposition}\label{P2}
Let $E_0,E_i\in\mathcal{C}^n$, $i\in\mathbb{N}$. If $E_i\rightarrow E_0$
in the $L_1$ distance and $P(E_i)\rightarrow P(E_0)$, then $\Pi E_i\rightarrow\Pi E_0$ in the
Hausdorff distance.
\end{proposition}
\noindent{\bf Proof.}
Since $E_i$ converges to $E_0$ in the $L_1$ distance when $i\rightarrow\infty$, we have
$$\chi_{E_i}\;{\rm converges\;to}\;\chi_{E_0}\;{\rm with\;respect\;to}\;L^1(\mathbb{R}^n).$$
 Since $|D\chi_{E_i}|(\mathbb{R}^n)=P(E_i)$ and $P(E_i)\rightarrow P(E_0)$, $|D\chi_{E_i}|(\mathbb{R}^n)$ is uniformly bounded. Hence, by \cite[Proposition 3.13]{MR1857292} one deduces that
$$\chi_{E_i}\rightharpoonup \chi_{E_0}\;\;{\rm weakly^{\ast}}\;{\rm in}\;\mathbb{R}^n\;{\rm when}\;i\rightarrow\infty.$$
Thus
\begin{eqnarray}\label{2y}
D\chi_{E_i}\rightharpoonup D\chi_{E_0}\;\;{\rm weakly^{\ast}}\;{\rm in}\;\mathbb{R}^n\;{\rm when}\;i\rightarrow\infty.
\end{eqnarray}

Since $P(E_i)\rightarrow P(E_0)$,
\begin{eqnarray}\label{2h}
|D\chi_{E_i}|(\mathbb{R}^n)\rightarrow |D\chi_{E_0}|(\mathbb{R}^n).
\end{eqnarray}

By (\ref{2y}), (\ref{2h}) and Reshetnyak continuity theorem (see \cite[Theorem 2.39]{MR1857292}), we have
\begin{eqnarray}\label{2i}
\lim_{i\rightarrow\infty}\int_{\mathbb{R}^n}\left|u\cdot\frac{D\chi_{E_i}(x)}{|D\chi_{E_i}(x)|}\right|d|D\chi_{E_i}|(x)=\int_{\mathbb{R}^n}\left|u\cdot\frac{D\chi_{E_0}(x)}{|D\chi_{E_0}(x)|}\right|d|D\chi_{E_0}|(x).
\end{eqnarray}

By (\ref{2b}), (\ref{2c}) and (\ref{2i}), for any $u\in S^{n-1}$, we have
\begin{eqnarray}\label{2j}
\lim_{i\rightarrow\infty}\int_{\partial^{\ast}E_i}\left|u\cdot\nu^{E_i}(x)\right|d\mathcal{H}^{n-1}(x)=\int_{\partial^{\ast}E_0}\left|u\cdot\nu^{E_0}(x)\right|d\mathcal{H}^{n-1}(x).
\end{eqnarray}

By the definition (\ref{3n}) of $\Pi E$ and (\ref{2j}), $\lim_{i\rightarrow\infty}h_{\Pi E_i}(u)=h_{\Pi E_0}(u)$.
Since the support functions $h_{\Pi E_i}\rightarrow h_{\Pi E_0}$ pointwise (on $S^{n-1}$) they converge unifromly (see, e.g., Schneider \cite[p. 54]{MR3155183}) completing the proof.\qed
\

We now demonstrate the affine nature of the Petty projection operator. Our proof follows along the same lines as that of Lemma 2.6 proved by Lutwak et al. \cite{MR2563216}.
\begin{proposition}\label{P1}
If $E\in\mathcal{C}^n$ and $A\in SL(n)$, then
$$\Pi^{\ast}AE=A\Pi^{\ast}E.$$
\end{proposition}
\noindent{\bf Proof.}
We first suppose that $P\in\mathcal{C}^n$ is a polyhedra satisfying
\begin{eqnarray}\label{2m}
\mathcal{H}^{n-1}(\{x\in\partial^{\ast}P:\;x\cdot\nu^{P}(x)=0\})=0.
\end{eqnarray}
Suppose $(n-1)$-dimensional faces of $P$ are $F_1,\dots,F_m$. Let $u_1,\dots,u_m$ be the outer unit normals to the faces, and let $h_1,\dots,h_m$ denote the {\it support numbers} of the faces of $P$; i.e., $h_i=|w_i\cdot u_i|$, where $w_i\in F_i$. Let $V_1,\dots,V_m$ denote the volumes of the facial cones, so that, $V_i=\frac{1}{n}|h_i||F_i|$. By (\ref{2m}), $h_i\neq 0$ for $i=1,\dots,m$.

For $A\in SL(n)$, let $P^{\diamond}=AP=\{Ax:\;x\in P\}$. Let $F_1^{\diamond},\dots,F_m^{\diamond}$ denote the faces of $P^{\diamond}$, let $u_1^{\diamond},\dots,u_m^{\diamond}$ be the outer unit normals of the faces of $P^{\diamond}$ and let $h_1^{\diamond},\dots,h_m^{\diamond}$ denote the corresponding support numbers of $P^{\diamond}$. Since $A\in SL(n)$, obviously the volumes $V_1^{\diamond},\dots,V_m^{\diamond}$ of the facial cones of $P^{\diamond}$ are such that $V_i^{\diamond}=V_i$.

The face $F_i$ parallel to the subspace $u_i^{\perp}$ is transformed by $A$
into the face $F_i^{\diamond}=A F_i$ parallel to $(A^{-t}u_i)^{\perp}$ and thus
\begin{eqnarray}\label{2g}
u_i^{\diamond}=A^{-t}u_i/\left|A^{-t}u_i\right|.
\end{eqnarray}
For $w_i^{\diamond}\in F_i^{\diamond}$, there exists $w_i\in F_i$ such that $w_i^{\diamond}=Aw_i$. Thus, from (\ref{2g}), we have
\begin{eqnarray}\label{2ff}
h_i^{\diamond}=|w_i^{\diamond}\cdot u_i^{\diamond}|=|Aw_i\cdot u_i^{\diamond}|=|w_i\cdot A^tu_i^{\diamond}|=|w_i\cdot u_i|/\left|A^{-t}u_i\right|=h_i/\left|A^{-t}u_i\right|.
\end{eqnarray}
Since $h_i\neq 0$, $h^{\diamond}_i\neq 0$ for $i=1,\dots,m$.

Now from definition (\ref{3n}), the fact that $V_i^{\diamond}=V_i$ together with (\ref{2g}) and (\ref{2ff}), definition (\ref{3n}) again, we have, for $z\in\mathbb{R}^n$,
\begin{eqnarray*}
h_{\Pi AP}(z)&=&h_{\Pi P^{\diamond}}(z)
=\frac{1}{2}\sum_{i=1}^{m}\left|z\cdot u_i^{\diamond}\right||F^{\diamond}_i|
=\frac{n}{2}\sum_{i=1}^{m}\left|\frac{z\cdot u_i^{\diamond}}{h^{\diamond}_i}\right||V^{\diamond}_i|\nonumber\\
&=&\frac{n}{2}\sum_{i=1}^{m}\left|\frac{z\cdot A^{-t}u_i}{h_i}\right||V_i|=\frac{n}{2}\sum_{i=1}^{m}\left|\frac{A^{-1}z\cdot u_i}{h_i}\right||V_i|\nonumber\\
&=&h_{\Pi P}(A^{-1}z)=h_{A^{-t}\Pi P}(z),
\end{eqnarray*}
showing that $\Pi AP=A^{-t}\Pi P$. By (\ref{1f}), $\Pi^{\ast}AP=A\Pi^{\ast}P$. This, together with Proposition \ref{P2} and Theorem \ref{T11}, completes the proof.\qed

\section{Proof of the main theorems}\label{S5}

\begin{lemma}\label{L1}\cite[Lemma 3.2]{MR2178968}
Let $E\in\mathcal{C}^n$. Then
\begin{eqnarray}\label{3ee}
\frac{\partial\ell}{\partial x_i}(x')=\int_{(\partial^{\ast}E)_{x'}}\frac{\nu_i^E(x',y)}{|\nu_y^E(x',y)|}d\mathcal{H}^0(y),\;\;\;i=1,\dots,n-1,
\end{eqnarray}
 for $\mathcal{L}^{n-1}$-a.e. $x'\in\pi_{n-1}(E)^+$.
\end{lemma}

\begin{remark}
An application of Lemma \ref{L1} and of (\ref{2ee}) to $E^s$ yields, in particular,
\begin{eqnarray}\label{3e}
\frac{\partial\ell}{\partial x_i}(x')=2\frac{\nu_i^{E^s}(x',\frac{1}{2}\ell(x'))}{\left|\nu_y^{E^s}(x',\frac{1}{2}\ell(x'))\right|}\;{\rm for}\;\mathcal{L}^{n-1}{\text-}{\rm a.e.}\;x'\in\pi_{n-1}(E)^+.
\end{eqnarray}
\end{remark}

\begin{lemma}\cite[Lemma 4.1]{MR3055761}\label{L10}
Let $E$ be any set of finite perimeter in $\mathbb{R}^n$, $n\geq 2$, and let $U$ be any Borel subset of $\mathbb{R}^{n-1}$. Then
$$\mathcal{H}^{n-1}(\{x\in\partial^{\ast} E:\;\nu_y^E(x)=0\}\cap(U\times\mathbb{R}_y))=0$$
if and only if
$$P(E;V\times\mathbb{R}_y)=0\;for\;each\;Borel\;subset\;V\;of\;U\;such\;that\;\mathcal{L}^{n-1}(V)=0.$$
\end{lemma}

\begin{lemma}\label{L4}
Let $E\in\mathcal{C}^n$ satisfy (\ref{1b}). Then
\begin{eqnarray}\label{2a}
\left(\Pi^{\ast} E\right)^s\subset\Pi^{\ast}E^s.
\end{eqnarray}
\end{lemma}
\noindent{\bf Proof.}
We will be appealing to Lemma \ref{L2} and thus we begin by supposing that
$$h_{\Pi E}(x_0',t)=1\;\;{\rm and}\;\;h_{\Pi E}(x_0',-s)=1,$$
with $t\neq -s$, or equivalently, by (\ref{2k}) and (\ref{2kk}), that
$$(x_0',t)\in\partial\Pi^{\ast} E\;\;{\rm and}\;\;(x_0',-s)\in\partial\Pi^{\ast} E.$$
By (\ref{3n})
\begin{eqnarray}\label{3c}
\frac{1}{2}\int_{\partial^{\ast}E} |(x_0',t)\cdot \nu^E(x)| d\mathcal{H}^{n-1}(x)=1
\end{eqnarray}
and
\begin{eqnarray}\label{3cc}
\frac{1}{2}\int_{\partial^{\ast}E} |(x_0',-s)\cdot \nu^E(x)| d\mathcal{H}^{n-1}(x)=1.
\end{eqnarray}
By Lemma \ref{L2}, the desired inclusion (\ref{2a}) will have been
established if we can show that
\begin{eqnarray}\label{3bb}
h_{\Pi E^s}\left(x_0',\frac{1}{2}t+\frac{1}{2}s\right)\leq1.
\end{eqnarray}

By \cite[Proposition 4.2]{MR2178968}, if (\ref{1b}) is established, then
\begin{eqnarray}\label{1bb}
\mathcal{H}^{n-1}(\{x\in\partial^{\ast}E^s:\;\nu_y^{E^s}(x)=0\})=0.
\end{eqnarray}

Let $G_E$ and $G_{E^s}$ be the sets associated with $E$ and $E^s$, respectively, as in Theorem \ref{T7}.
By $G_E,G_{E^s}\subset\pi_{n-1}(E)^+$, $\mathcal{L}^{n-1}(\pi_{n-1}(E)^+\setminus G_{E^s})=0$, $\mathcal{L}^{n-1}(\pi_{n-1}(E)^+\setminus G_{E})=0$, (\ref{1b}), (\ref{1bb}), Lemma \ref{L10}, (\ref{3n}), (\ref{3c}) and (\ref{3cc}), we have that (\ref{3bb}) is equivalent to
\begin{eqnarray}\label{3bbb}
&&\int_{\partial^{\ast}E^s\cap (G_{E^s}\times\mathbb{R})} \left|\left(x_0',\frac{t+s}{2}\right)\cdot \nu^{E^s}(x)\right| d\mathcal{H}^{n-1}(x)\\
&&\leq\frac{1}{2}\int_{\partial^{\ast}E\cap (G_E\times\mathbb{R})} |(x_0',t)\cdot \nu^E(x)| d\mathcal{H}^{n-1}(x)+\frac{1}{2}\int_{\partial^{\ast}E\cap (G_E\times\mathbb{R})} |(x_0',-s)\cdot \nu^E(x)| d\mathcal{H}^{n-1}(x)\nonumber.
\end{eqnarray}

We have the following chain of equalities:
\begin{eqnarray}\label{3k}
&&\int_{\partial^{\ast}E^s\cap(G_{E^s}\times\mathbb{R})} \left|\left(x_0',\frac{t+s}{2}\right)\cdot \nu^{E^s}(x)\right| d\mathcal{H}^{n-1}(x)\\
&=&\int_{G_{E^s}} dx'\int_{(\partial^{\ast}E^s)_{x'}} \left|\left(x_0',\frac{t+s}{2}\right)\cdot \nu^{E^s}(x)\right| \frac{d\mathcal{H}^0(y)}{\left|\nu^{E^s}_y(x',y)\right|}\nonumber\\
&=&\int_{G_E} dx'\int_{(\partial^{\ast}E^s)_{x'}} \left|\left(x_0',\frac{t+s}{2}\right)\cdot \left(\frac{\nu_{x'}^{E^s}(x)}{{\left|\nu^{E^s}_y(x',y)\right|}},\frac{\nu_y^{E^s}(x)}{{\left|\nu^{E^s}_y(x',y)\right|}}\right)\right| d\mathcal{H}^0(y),\nonumber
\end{eqnarray}
where the first is due to the co-area formula (\ref{2l}) and (\ref{2gg}), the second to the fact that $\mathcal{L}^{n-1}(G_{E^s}\triangle G_E)=0$.

By (\ref{3e}) and (\ref{2ee}), for $E^s$
\begin{eqnarray}\label{3g}
&&\int_{G_E} dx'\int_{(\partial^{\ast}E^s)_{x'}} \left|\left(x_0',\frac{t+s}{2}\right)\cdot \left(\frac{\nu_{x'}^{E^s}(x)}{{\left|\nu^{E^s}_y(x',y)\right|}},\frac{\nu_y^{E^s}(x)}{{\left|\nu^{E^s}_y(x',y)\right|}}\right)\right| d\mathcal{H}^0(y)\\
&=&\int_{G_E} dx'\int_{\partial^{\ast}(E^s)_{x'}} \left|\left(x_0',\frac{t+s}{2}\right)\cdot \left(\frac{1}{2}\nabla\ell(x'),\frac{\nu_y^{E^s}(x)}{{\left|\nu^{E^s}_y(x',y)\right|}}\right)\right| d\mathcal{H}^0(y)\nonumber\\
&=&\int_{G_E} dx'\int_{\partial^{\ast}(E^s)_{x'}\cap (G_E\times\mathbb{R}^+)} \left|\left(x_0',\frac{t+s}{2}\right)\cdot \left(\frac{1}{2}\nabla\ell(x'),-1\right)\right| d\mathcal{H}^0(y)\nonumber\\
&&+\int_{G_E} dx'\int_{\partial^{\ast}(E^s)_{x'}\cap (G_E\times\mathbb{R}^-)} \left|\left(x_0',\frac{t+s}{2}\right)\cdot \left(\frac{1}{2}\nabla\ell(x'),1\right)\right| d\mathcal{H}^0(y)\nonumber\\
&=&\frac{1}{2}\int_{G_E}\left|\left(x_0',t+s\right)\cdot \left(\nabla\ell(x'),-1\right)\right| dx'+\frac{1}{2}\int_{G_E} \left|\left(x_0',t+s\right)\cdot \left(\nabla\ell(x'),1\right)\right| dx'. \nonumber
\end{eqnarray}
For $x'\in G_E$, let
\begin{eqnarray}\label{3j}
\partial^{r,\ast}E_{x^{\prime}}:=\left\{(x^{\prime},y)\in\partial^{\ast}E_{x^{\prime}}:\;\nu^E_y(x^{\prime},y)<0\right\}
\end{eqnarray}
and
\begin{eqnarray}\label{3jj}
\partial^{l,\ast}E_{x^{\prime}}:=\left\{(x^{\prime},y)\in\partial^{\ast}E_{x^{\prime}}:\;\nu^E_y(x^{\prime},y)>0\right\}.
\end{eqnarray}
By Remark \ref{3s},
\begin{eqnarray}\label{3f}
\mathcal{H}^0(\partial^{l,\ast}(E_{x^{\prime}}))=\mathcal{H}^0(\partial^{r,\ast}(E_{x^{\prime}}))=m(x^{\prime})=\frac{\mathcal{H}^0(\partial^{\ast}(E_{x^{\prime}}))}{2}.
\end{eqnarray}
For $k=1,\dots,m(x^{\prime})$, let
 \begin{eqnarray}
 y_{k,l}(x^{\prime})\;{\rm be\; the\;}k{\text-}{\rm th}\;{\rm number}\;y{\rm\;satisfying}\;(x^{\prime},y)\in \partial^{l,\ast}E_{x^{\prime}},
 \end{eqnarray}
 \begin{eqnarray}\label{3ff}
 y_{k,r}(x^{\prime})\;{\rm be\; the\;}k{\text-}{\rm th}\;{\rm number}\;y{\rm\;satisfying}\;(x^{\prime},y)\in \partial^{r,\ast}E_{x^{\prime}}.
 \end{eqnarray}

By (\ref{3ee}), the monotonicity of $|a+br|+|a-br|$ with respect to $r>0$, where $a,b\in\mathbb{R}$, (\ref{3f})-(\ref{3ff}), the last expression of (\ref{3g})
\begin{eqnarray}\label{3h}
&\leq&\frac{1}{2}\int_{G_E}\left|(x_0',t+s)\cdot\left(\int_{\partial^{\ast}E_{x^{\prime}}}\frac{\nu_{x'}^E}{|\nu^E_y|}d\mathcal{H}^0,-\frac{1}{2}\int_{\partial^
{\ast}E_{x^{\prime}}}d\mathcal{H}^0\right)\right|dx^{\prime}\\
&&+\frac{1}{2}\int_{G_E}\left|(x_0',t+s)\cdot\left(\int_{\partial^{\ast}E_{x^{\prime}}}\frac{\nu_{x'}^E}{|\nu^E_y|}d\mathcal{H}^0,\frac{1}{2}\int_{\partial^
{\ast}E_{x^{\prime}}}d\mathcal{H}^0\right)\right|dx^{\prime}\nonumber\\
&=&\frac{1}{2}\int_{G_E}\left|(x_0',t+s)\cdot\left(\sum^{m(x^{\prime})}_{k=1}\left(\frac{\nu_{x'}^E}{|\nu^E_y|}(x^{\prime},y_{k,l}(x^{\prime}))+\frac{\nu_{x'}^E}{|\nu^E_y|}(x^{\prime},y_{k,r}(x^{\prime}))\right),-m(x^{\prime})\right)\right|dx^{\prime}\nonumber\\
&&+\frac{1}{2}\int_{G_E}\left|(x_0',t+s)\cdot\left(\sum^{m(x^{\prime})}_{k=1}\left(\frac{\nu_{x'}^E}{|\nu^E_y|}(x^{\prime},y_{k,l}(x^{\prime}))+\frac{\nu_{x'}^E}{|\nu^E_y|}(x^{\prime},y_{k,r}(x^{\prime}))\right),m(x^{\prime})\right)\right|dx^{\prime}.\nonumber
\end{eqnarray}

By the convexity of the absolute value function $|\cdot|$, (\ref{3j}), (\ref{3jj}) and the co-area formula (\ref{2l}), the last expression of (\ref{3h})
\begin{eqnarray}\label{3i}
&\leq&\frac{1}{2}\int_{G_E}\sum_{k=1}^{m(x^{\prime})}\left|(x_0',t)\cdot\left(\frac{\nu_{x'}^E}{|\nu^E_y|}(x^{\prime},y_{k,r}(x^{\prime})),-1\right)\right|dx^{\prime}\\
&&+\frac{1}{2}\int_{G_E}\sum_{k=1}^{m(x^{\prime})}\left|(x_0',-s)\cdot\left(\frac{\nu_{x'}^E}{|\nu^E_y|}(x^{\prime},y_{k,l}(x^{\prime})),1\right)\right|dx^{\prime}\nonumber\\
&&+\frac{1}{2}\int_{G_E}\sum_{k=1}^{m(x^{\prime})}\left|(x_0',t)\cdot\left(\frac{\nu_{x'}^E}{|\nu^E_y|}(x^{\prime},y_{k,l}(x^{\prime})),1\right)\right|dx^{\prime}\nonumber\\
&&+\frac{1}{2}\int_{G_E}\sum_{k=1}^{m(x^{\prime})}\left|(x_0',-s)\cdot\left(\frac{\nu_{x'}^E}{|\nu^E_y|}(x^{\prime},y_{k,r}(x^{\prime})),-1\right)\right|dx^{\prime}\nonumber\\
&=&\frac{1}{2}\int_{G_E}dx^{\prime}\int_{\partial^{\ast}E_{x^{\prime}}}\left|(x_0',t)\cdot\left(\frac{\nu_{x'}^E}{|\nu^E_y|},\frac{\nu^E_y}{|\nu^E_y|}\right)\right|d\mathcal{H}^0(y)\nonumber\\
&&+\frac{1}{2}\int_{G_E}dx^{\prime}\int_{\partial^{\ast}E_{x^{\prime}}}\left|(x_0',-s)\cdot\left(\frac{\nu_{x'}^E}{|\nu^E_y|},\frac{\nu^E_y}{|\nu^E_y|}\right)\right|d\mathcal{H}^0(y)\nonumber\\
&=&\frac{1}{2}\int_{\partial^{\ast}E}\left|(x_0',t)\cdot\nu^{E}\right|d\mathcal{H}^{n-1}+\frac{1}{2}\int_{\partial^{\ast}E}\left|(x_0',-s)\cdot\nu^{E}\right|d\mathcal{H}^{n-1}\nonumber.
\end{eqnarray}
From (\ref{3k}), (\ref{3g}), (\ref{3h}) and (\ref{3i}) we get (\ref{3bbb}).  This completes the proof of the lemma.\qed
\

\noindent{\bf Proof of Theorem \ref{T1}.} By Lemma \ref{L4} and the volume invariance of Steiner symmetrization, inequality (\ref{1a}) is established. \qed
\

Next, we prove Theorem \ref{T2}. The following lemma was proved in  \cite[Lemma 3.1]{MR3085623} (also see \cite{MR367161}). We give a different proof here.
\begin{lemma}\label{L3}
 Let $E,E_i\in \mathcal{C}^n$, $i\in\mathbb{N}$. If $|E|=|E_i|$ and $E_i\rightarrow E$ in the Hausdorff distance, then
$E_i\rightarrow E$ in the $L_1$ distance.
\end{lemma}
\noindent{\bf Proof.}
Since $E_i\rightarrow E$ in the Hausdorff distance, for any $\varepsilon>0$, there exists a positive integer $N$ such that for $i>N$,
$$E_i\subset E+\varepsilon B^n.$$
Thus
\begin{eqnarray}\label{5a}
E_i\backslash E\subset (E+\varepsilon B^n)\backslash E.
\end{eqnarray}

 For a decreasing sequence of positive numbers $\varepsilon_k\rightarrow 0^+$, by the limit theorem with respect to  sequences of measurable sets (see \cite[Theorem 1.2 (iv)]{MR3409135}), we have
\begin{eqnarray}\label{5aa}
\lim\limits_{k\rightarrow\infty}\left|(E+\varepsilon_kB^n)\backslash E\right|=\left|\bigcap_{k=1}^{\infty}(E+\varepsilon_k B^n)\backslash E\right|=0.
\end{eqnarray}
By (\ref{5a}) and (\ref{5aa}),
\begin{eqnarray}\label{5b}
\lim_{i\rightarrow \infty}|E_i\backslash E|=0.
\end{eqnarray}

Since $|E|=|E_i|$, we have
\begin{eqnarray}\label{5c}
|E\backslash E_i|=|E|-|E\cap E_i|=|E_i|-|E\cap E_i|=|E_i\backslash E|.
\end{eqnarray}
By (\ref{5b}) and (\ref{5c}), $\lim_{i\rightarrow \infty}d_1(E_i,E)=2|E_i\backslash E|=0$.\qed
\

The Petty projection operator $\Pi:\mathcal{C}^n\rightarrow\mathcal{K}_0^n$ is
weakly continuous in some sense when $E_i$ converges to $E^{\star}$ in the $L_1$ distance without the assumption that $P(E_i)\rightarrow P(E^{\star})$.

\begin{lemma}\label{L5}
Let $E,E_i\in\mathcal{C}^n$, $i\in\mathbb{N}$. If $L_1(E_i,E^{\star})\rightarrow 0$,
then there exist a subsequence of $\{\Pi E_i\}_{i=1}^{\infty}$, denoted by $\{\Pi E_i\}_{i=1}^{\infty}$ as well, and a convex body $K$  such that $0\in K$,
$d_H\left(\Pi E_i,K\right)\rightarrow 0$
and $\Pi E^{\star}\subset K$.
\end{lemma}

\noindent{\bf Proof.} By the definition (\ref{3n}) of $\Pi E$, (\ref{2n}) and (\ref{2c}), for any $u\in S^{n-1}$,
$$h_{\Pi E_i}(u)=\frac{1}{2}\int_{\partial^{\ast}E_i} |u\cdot \nu^{E_i}(x)| d\mathcal{H}^{n-1}(x)\leq\frac{1}{2}\mathcal{H}^{n-1}(\partial^{\ast}E_i)=\frac{1}{2}P(E_i).$$
Since $P(E_i)$ is decreasing with respect to $i$ (see Theorem \ref{T9}), there exists a constant $r_0>0$ such that
$\Pi E_i\subset r_0B^n\;\;{\rm for\;any}\;i$.
 By Blaschke selection theorem (see \cite[Theorem 1.8.7]{MR3155183}), there exists a subsequence of $\{\Pi E_i\}_{i=1}^{\infty}$, denoted by $\{\Pi E_i\}_{i=1}^{\infty}$ as well, that converges to a convex body $K$ in the Hausdorff distance. Since $h_{\Pi E_i}(u)>0$ for any $i$ and $u\in S^{n-1}$ and $h_{\Pi E_i}(u)\rightarrow h_K(u)$, $h_K\geq0$. Thus, $0\in K$.

Since $\{E_i\}_{i=1}^{\infty}$ converges to $E^{\star}$ in the $L_1$ distance when $i\rightarrow\infty$, we have
$\chi_{E_i}$ converges to $\chi_{E^{\star}}$ with respect to $L_1(\mathbb{R}^n)$.
 Since $|D\chi_{E_i}|(\mathbb{R}^n)=P(E_i)$ and $P(E_i)$ is decreasing with respect to $i$ (see Theorem \ref{T9}), $|D\chi_{E_i}|(\mathbb{R}^n)$ is uniformly bounded. Hence, by \cite[Proposition 3.13]{MR1857292} one deduces that
\begin{eqnarray}\label{6e}
D\chi_{E_i}\rightharpoonup D\chi_{E^{\star}}\;\;{\rm weakly^{\ast}}\;{\rm in}\;\mathbb{R}^n\;{\rm when}\;i\rightarrow\infty.
\end{eqnarray}

By (\ref{6e}) and Reshetnyak lower semicontinuity theorem (see \cite[Theorem 2.38]{MR1857292}), we have
\begin{eqnarray}\label{6f}
\int_{\mathbb{R}^n}\left|u\cdot\frac{D\chi_{E^{\star}}(x)}{|D\chi_{E^{\star}}(x)|}\right|d|D\chi_{E^{\star}}|(x)\leq\lim_{i\rightarrow\infty}\int_{\mathbb{R}^n}\left|u\cdot\frac{D\chi_{E_i}(x)}{|D\chi_{E_i}(x)|}\right|d|D\chi_{E_i}|(x).
\end{eqnarray}

By (\ref{2b}), (\ref{2c}) and (\ref{6f}), for any $u\in S^{n-1}$, we have
\begin{eqnarray}\label{6g}
\int_{\partial^{\ast}E^{\star}}\left|u\cdot\nu^{E^{\star}}(x)\right|d\mathcal{H}^{n-1}(x)\leq\lim_{i\rightarrow\infty}\int_{\partial^{\ast}E_i}\left|u\cdot\nu^{E_i}(x)\right|d\mathcal{H}^{n-1}(x).
\end{eqnarray}

By the definition (\ref{3n}) of $\Pi E$ and (\ref{6g}), $h_{\Pi E^{\star}}(u)\leq\lim_{i\rightarrow\infty}h_{\Pi E_i}(u)=h_K(u)$. Thus, $\Pi E^{\star}\subset K$.\qed
\

\noindent{\bf Proof of Theorem \ref{T2}.} By the monotonicity of the projection operator (Theorem \ref{T1}), the convergence of the subsequence $\{D_k\}$ of $\{E_j\}$ (Theorem \ref{T3} ), the weakly continuity of the projection operator (Lemma \ref{L5}), for $E\in\mathcal{C}^n$, the
volume $V(E)^{n-1}V(\Pi^{\ast}E)$ is maximized when $E=E^{\star}$. Since $V(E^{\star})^{n-1}V(\Pi^{\ast}E^{\star})=(\omega_n/\omega_{n-1})^n$, inequality (\ref{1i}) is established. \qed

\section{Open problems}

\begin{problem}\label{Pr1}
For $E\in\mathcal{C}^n$, if $E_i:=S_{u_i}\dots S_{u_1}E$ for a sequence of directions $\{u_i\}$ and there exists a subsequence $E_{i_j}$ converges to $E^{\star}$ in the Hausdorff distance, then does $E_i$ converge to $E^{\star}$ in the Hausdorff distance?
\end{problem}

If the answer of Problem \ref{Pr1} is positive, then the sequence $\{E_i\}$ in Theorem \ref{T3} converges to $E^{\star}$, which is stronger than the convergence of its subsequence.

\begin{problem}\label{Pr2}
For $E\in\mathcal{C}^n$, if $E_i:=S_{u_i}\dots S_{u_1}E$ for a sequence of directions $\{u_i\}$ and $E_i\rightarrow E^{\star}$ in the Hausdorff distance, then does $P(E_i)$ converge to $P(E^{\star})$?
\end{problem}

If the answer of Problem \ref{Pr2} is positive, then by Proposition \ref{P2}, the sequence $\{\Pi E_i\}$ in Lemma \ref{L5} converges to $\Pi E^{\star}$, which is stronger than its convergence to a convex body $K$ satisfying $\Pi E^{\star}\subset K$.



%
%


\bibliographystyle{plain}

\end{document}